\documentclass[12pt]{article}
\usepackage{epsfig}
\font\tc=wncyr10

\input cyracc.def
\usepackage{amsmath}
\usepackage{amsfonts}
\usepackage{latexsym}
\usepackage{graphicx, epsfig}

\def\dist{{\mathrm{dist}\,}}
\def\clos{{\mathrm{clos}\,}}
\def\U{{\mathbf{U}}}
\def\D{{\mathcal{D}}}
\def\es{{\mathcal{S}}}
\title{Oscillation of Fourier Integrals with a spectral gap}
\author{A. Eremenko\thanks{Supported by NSF grant DMS 0100512 and
by the Humboldt Foundation.}$\;$ and D. Novikov\thanks{Supported
by NSF grant DMS 0200861 and by the McDonnell Foundation.}}
\date{January 2, 2003}
\newtheorem{theorem}{Theorem}
\newtheorem{lemma}{Lemma}
\newtheorem{proposition}{Proposition}
\newtheorem{statement}{Statement}

\newtheorem{example}{Example}

\def\E{\mathbf{E}}

\def\T{{\mathbf{T}}}
\def\H{{\mathcal{H}}}
\def\C{\mathbf{C}}
\def\bC{{\overline{\mathbf C}}}
\def\R{{\mathbf{R}}}
\def\Z{{\mathbf{Z}}}
\def\F{{\mathcal{F}}}
\def\Fr{{\mathrm{Fr}}}
\begin{document}
\maketitle

\section{Introduction}

Suppose that in a real Fourier series, first $m$ terms vanish:
\begin{equation}
\label{poly}
f(x)=\sum_{n\geq m}(c_ne^{inx}+\overline{c}_ne^{-inx}),\quad f\neq 0.
\end{equation}
Then $f$ has at least $2m$ changes of sign on the interval $|x|\leq\pi$.
For trigonometric polynomials this follows from
a result of Sturm \cite{Sturm}; the general case is due to
Hurwitz.

Here is a simple proof. The number
of sign changes is even. If $f$ has at most $2(m-1)$ changes
of sign then we can find a trigonometric polynomial $g$ of degree at most $m-1$
which changes sign at the same places as $f$. Then $fg$ is
of constant sign
which contradicts the orthogonality of $f$ and $g$.

We consider the following extension of this result to Fourier integrals.
\begin{statement}\label{prop0}
Suppose that a real function $f$ has a spectral gap, that is
its Fourier transform is zero on an
interval $(-a,a)$. Then the number of sign changes $s(r,f)$ of $f$
on the interval $[0,r]$ satisfies
\begin{equation}\label{lowdens}
\liminf_{r\to\infty}\frac{s(r,f)}{r}\geq\frac{a}{\pi}.
\end{equation}
\end{statement}
The proof given above permits to estimate the {\em upper} density
of the sequence of sign changes of $f$, but our goal is to estimate
the {\em lower} density.

In our discussion of Statement 1 it is convenient to use the
general definition of Fourier transform which is due to Carleman
\cite{CarlemanF}. Suppose that a locally integrable function $f$
satisfies
\begin{equation}
\label{exp}
\int |f(x)|e^{-\lambda |x|}dx<\infty\quad\mbox{for all}\quad\lambda>0.
\end{equation}
Then the functions
$$F^+(z)=\int_{-\infty}^0f(t)e^{-itz}dt\quad\mbox{and}\quad
F^-(z)=-\int_{0}^\infty f(t)e^{-itz}dt$$
are analytic in the upper and lower half-planes, respectively,
and the {\em generalized Fourier transform} is defined as a
{\em hyperfunction}, that is a
pair
$(F^+,F^-)$.
modulo addition of an entire function.
Fourier transform in the usual sense, if exists, is obtained
as the difference of the boundary values,
$\hat{f}(x)=F^+(x)-F^-(x),\; x\in\R$. If $f$ is a temperate distribution,
then the boundary values in the previous formula can be
interpreted as limits of temperate distributions.
Thus suggests a general definition of the spectrum of
a function.
\vspace{.1in}

\noindent
{\bf Definition 1}. {\em The spectrum of a function $f$ satisfying $(\ref{exp})$
is
the complement of the maximal open set $U\subset \R$ such that
$F^+$ and $F^-$ are analytic continuations of each other
through $U$.}
\vspace{.1in}

Thus a function $f$ has a spectral gap $(-a,a)$ if $F^+$ and $F^-$ are
analytic continuations of each other through the interval
$(-a,a)$.

In engineering
literature, functions with a spectral gap are called high-pass signals.

Condition (\ref{exp}) is too weak to develop a proper generalization
of Harmonic Analysis \cite{Beur3}, for example the spectrum of
a function satisfying (\ref{exp}) can be empty.
Of the many generalizations
of classical theories of Fourier transform
we mention first of all the theory of
temperate distributions of Schwartz \cite{Hormander}.
A further generalization was proposed
by Beurling in his lectures \cite{Beurling}\footnote{\noindent
Unfortunately,
\cite{Beurling}
is unpublished. There is an exposition of
Beurling's theory in~\cite{Bjorck}.}.

Following Beurling
\cite{Beur3,Beurling}, we consider locally integrable functions
$f$ which satisfy
\begin{equation}
\label{beu}
\int|f(x)|e^{-\lambda\omega(x)}dx<\infty\quad\mbox{for some}\quad \lambda>0,
\end{equation}
where $\omega\geq 1$ is a real function with the property
\begin{equation}
\label{weight}
\int\frac{\omega(x)}{1+x^2}dx<\infty.
\end{equation}
Suppose, in addition, that $\omega(x)\geq \log(1+|x|),\; x\in\R$, and
\begin{equation}\label{weight1}
0=\omega(0)\leq \omega(x+y)\leq\omega(x)+\omega(y),\quad x,y\in\R.
\end{equation}
The space $\es_\omega$
of {\em test functions} consists of all functions $\phi$ in $L^1:=L^1(\R)$,
such that $\phi$ and its Fourier transform
$$\hat{\phi}(t)=\int\phi(x)e^{-itx}dx$$
belong
to $C^\infty$ and satisfy
\begin{equation}\label{seminorms1}
\sup_{\R}|\phi^{(k)}|e^{\lambda\omega}<\infty,\quad
\sup_\R|\hat{\phi}^{(k)}|e^{\lambda\omega}<\infty,
\end{equation}
for all non-negative integers $k$ and all $\lambda\geq 0$. The topology
on $\es_\omega$ is defined by the seminorms
(\ref{seminorms1}).
The dual space $\es_\omega^\prime$ is called the space of
$\omega$-{\em temperate distributions}. When $\omega=\log(1+|x|)$
we obtain the space $\es'$ of Schwartz's temperate distributions.
Fourier transform of a distribution $f$ is defined by
$$(\hat{f},\phi)=(f,\hat{\phi}).$$
The support of a distribution $f$ is the complement of the maximal open
set $U\subset\R$ such that $(f,\phi)=0$ for all $\phi\in\es_\omega$
with support in
$U$.
A complex-valued locally integrable function $f$ on the real line defines
an
$\omega$-tempered distribution if it satisfies (\ref{beu}).
Such functions $f$ also satisfy (\ref{exp}), and
new definition of support of $\hat{f}$
is consistent with more general Definition 1.

Conditions (\ref{beu}) and (\ref{weight}) imply
\begin{equation}
\label{logint}
\int\frac{\log^+|f(x)|}{1+x^2}dx<\infty.
\end{equation}
This property or, more precisely, the property
(\ref{weight}) of the weight $\omega$ ensures that test functions in
Beurling's theory {\em are not} quasianalytic, in particular,
there exist test functions with bounded support.
Functions with bounded spectrum which satisfy
(\ref{logint}) form a subclass of efet called the {\em Cartwright class}.

We will use weaker regularity assumptions about $\omega$ than (\ref{weight1}).
A real function $\omega\geq 1$ on the real line is called a
{\em Beurling--Malliavin weight} (BMW) if it satisfies (\ref{weight}) and,
in addition, has at least one of the
following properties:
\vspace{.1in}

\noindent
(i) $\omega$ is uniformly continuous, or
\vspace{.1in}

\noindent
(ii) $\exp\omega$ is the restriction of an entire function of
exponential type to the real line.
\vspace{.1in}

Notice that (\ref{weight1}) implies (i). Moreover, for a BMW $\omega$,
(\ref{beu}) implies (\ref{exp}).

Our main result shows that Statement 1 holds for $\omega$-temperate
distributions:
\begin{theorem}\label{thm0}
Let $\omega$ be a BMW.
If
$f\neq 0$ is a real measurable function satisfying $(\ref{beu})$
and having a spectral gap $(-a,a)$,
then
$(\ref{lowdens})$ holds.
\end{theorem}
In particular, this applies to locally integrable temperate distributions
of the Schwartz space $\es'$,
which contains,
for example,
all bounded functions.

The theory of mean motion \cite{Jessen,Levin} suggests a
stronger version of (\ref{lowdens}):
\begin{equation}\label{uniform}
\liminf_{x-y\to+\infty}\frac{n(x,f)-n(y,f)}{x-y}\geq \frac{a}{\pi}.
\end{equation}
This is not true, even for bounded $L^1$ functions with bounded spectrum:
\begin{example}\label{thm2} For every pair of positive numbers $a<b$,
there exists a real entire function $f$ of
exponential type $b$, whose restriction to the
real line is bounded and belongs to $L^1$, which has
a spectral gap $(-a,a)$, and the property that for a sequence of
intervals $[y_k,x_k]$ whose lengths tend to infinity,
$f$ has no zeros on $[y_k,x_k]$.
\end{example}

Examples of functions with a spectral gap and no sign changes
on {\em one} long interval
are contained in \cite{Logan}.

To show that condition (\ref{beu}) is essential in Theorem 1,
we consider functions $f$ with {\em bounded spectrum}.
This means that the generalized Fourier transform $(F_1,F_2)$
extends to a function $F$ analytic in $\bC\backslash[-b,b]$,
where $b\geq 0$. A theorem of P\'olya \cite{Bieberbach,Levin} gives a precise
description of such functions $f$: they are restrictions on the real
line of entire functions of exponential type $b$.
In engineering literature such functions are called band-limited
signals. We abbreviate ``entire function of exponential type'' as
{\em efet}.

\begin{example}\label{thm4}
For every positive numbers $a<b$, there exists
a real efet $f$ satisfying $(\ref{exp})$ whose spectrum
is
contained in $[-b,-a]\cup[a,b]$, and such that
$$\liminf_{r\to\infty} s(r,f)/r<a/\pi.$$
\end{example}

We conclude that condition (\ref{logint}) is essential for
validity of Statement 1.
Convergence or divergence of the integral (\ref{logint}) is a
fundamental dichotomy in Harmonic Analysis,
\cite{Beurling,Koosis}. From our point of view, the main
difference between the functions that satisfy (\ref{logint}) and
those that do not is explained by a theorem of Cartwright and Levinson
\cite[Ch.
5, Thm. 7]{Levin}: condition (\ref{logint}) implies completely regular
growth in the sense of Levin and Pfluger.
The situation is somewhat similar to the failure
of Titchmarsh's theorem
on the support of convolution \cite{R,KR} in the absence of condition
(\ref{logint}).

It is easy to construct examples of bounded
functions $f$ with bounded spectrum,
for which
the limit in (\ref{lowdens})
does not exist.
However, the theory of mean motion suggests the following question:
{\em under what additional conditions
does the limit in $(\ref{lowdens})$ exist? Does it exist
for trigonometric sums}
$$f(x)=\sum_{n=0}^m a_n\cos{\lambda_nx}+b_n\sin{\lambda_nx},
\quad \lambda_n, a_n, b_n\in\R
\;\mbox{?}$$

The paper is organized as follows. In section 2 we discuss known
results and conjectures about oscillation of functions with a spectral gap,
in section 3 we reduce our Theorem~1 to its special case
that $f\in L^1$, and construct Example~1.
In section 4 we prove Theorem 1 under the additional
assumption that $f$ is real analytic and has only simple zeros
on the real line. The general case is deduced in sections
5--7 by a smoothing procedure. Sections 8 and 9 are independent
of the rest of the paper. In section 8 we give a brief account
of Azarin's generalization of the theory of completely regular growth,
which we need for construction Example 2 in Section 9.

We thank Andrei Gabrielov, Iosif Ostrovskii, Misha Sodin,
and Serge Tabachnikov for valuable discussions, and
Jane Kinkus for procuring a copy of \cite{Beurling} for us.
The first-named author thanks
Tel-Aviv University where this work was completed.

\section{History and related results}

High-pass signals are important in Electrical Engineering.
Statement~\ref{prop0} was conjectured by Logan in his 1965 thesis
\cite{Logan} where he proved (\ref{lowdens}) under the additional
assumption that $f$ has bounded spectrum and is bounded on the
real line. One can replace in his result the condition of
boundness on the real line by the weaker condition (\ref{logint}).
So we have the following special case of our Theorem 1.
\begin{proposition}\label{logans} (Logan) Let $f\neq 0$ be a real
function
with bounded spectrum, having a spectral gap $(-a,a)$
and satisfying $(\ref{logint})$. Then $(\ref{lowdens})$ holds.
\end{proposition}
Example 2 shows that condition (\ref{logint}) cannot be dropped.
We include a proof for three reasons:
first, it is simple and gives a new proof
of Sturm's theorem itself,
second, we relax Logan's assumptions,
and third, his thesis is not everywhere easily available.

{\em Proof}. Let $b$ be the exponential type (bandwidth) of $f$, $b\geq a$.
As $f$ is real, it can be written as a sum
\begin{equation}\label{sum}
f(x)=h(x)+\overline{h}(x),\quad\mbox{where}\quad \overline{h}(z):=
\overline{h(\overline{z})},
\end{equation}
and $h$ is a function with a spectrum on $[a,b]$, which
satisfies (\ref{logint}). For the proof of this representation (\ref{sum})
see Proposition~\ref{111} in the next section. Now
\begin{equation}\label{logn}
f=e^{ibx}h_1+e^{-ibx}{\overline{h_1}}=\cos(bx)(h_1+{\overline{h_1}})
+i\sin(bx)(h_1-{\overline{h_1}}),
\end{equation}
where $h_1$ and ${\overline{h_1}}$ have their spectra on $[a-b,0]$ and $[0,b-a]$,
respectively. We conclude that $g=h_1+{\overline{h_1}}$
is a real efet with spectrum on
the interval $[(a-b),(b-a)]$, and $g$ satisfies (\ref{logint}).
Thus
by the theorem of Cartwright and Levinson \cite[Ch. V, Thm. 7]{Levin}, $g$ is an efet
of completely regular growth
in the sense of Levin and Pfluger.
In particular, the sequence of complex zeros of $g$ in any open
angle containing the positive ray has a density
equal to $(b-a)/\pi$. So the upper density of
positive zeros of $g$ is at most $(b-a)/\pi$. On the other hand, (\ref{logn})
implies
$$f(n\pi/b)=(-1)^ng(n\pi/b),$$
from which it is easy to derive that $s(r,f)\geq [br/\pi]-s(r,g)$. Dividing
by $r$ and passing to the lower limit, we obtain (\ref{lowdens}).
\hfill$\Box$\vspace{.1in}

It is important for this proof that $f$ is an efet.
Our Theorem~1, whose proof is based on different ideas,
extends Logan's result to functions with unbounded
spectrum.

The following conjecture of P.G.~Gri\-ne\-vich is contained in
\cite[(1996-5)]{Arnold}:
{\em ``If a real Fourier integral $f$ has a
spectral gap $(-a,a)$
then the limit average density of zeros of $f$ is at least $a/\pi$''.}

In the commentary to this problem in \cite{Arnold}, S.B.~Kuksin
mentioned the following result as a supporting evidence for
Grinevich's conjecture. Let $\xi(t)$ be a Gaussian stationary
random process, normalized by $\E\xi(0)=0$ and $\E\xi(0)^2=1$,
where $\E$ stands for the expectation. Let $r$ be the correlation function
of this process, $r(t)=\E\xi(0)\xi(t)$. Assume that the function
$r$ is integrable and has a spectral gap $(-a,a)$. Denote by
\def\e{{\mathcal{E}}}
$\e_T$ the random variable which is equal to the number of zeros of
the random function $\xi(t)$ on $[0,T]$. Then almost surely
$T^{-1}\e_T$ has a limit as $T\to\infty$, and this limit is at
least~$a/\pi$.

Other known results deal with averaged densities,
like
$$S(r,f)=\int_0^r\frac{s(t,f)+s(-t,f)}{t}dt.$$
When one uses $S(r,f)$,
condition (\ref{logint}) apparently plays no role anymore.
To demonstrate this, we
state and prove a version of Proposition \ref{logans}:
\begin{proposition}\label{logans2}
Let $f\neq 0$ be a real function with
bounded spectrum and a spectral gap
$(-a,a)$. Then
\begin{equation}
\label{log2}
\liminf_{r\to\infty}\frac{S(r,f)}{r}\geq \frac{2a}{\pi}.
\end{equation}
\end{proposition}
This property is weaker than (\ref{lowdens}).

{\em Proof}. We repeat the proof of Proposition \ref{logans}, but
instead of using the theorem of Cartwright and Levinson, apply
Jensen's formula. Decomposition (\ref{logn}) still holds, and
$g=h_1+ {\overline{h_1}}$ is an efet with the spectrum on
$[(a-b),(b-a)]$.
 Let $n(t,g)$ be the number of zeros of $g$
in the disc $\{ z:|z|\leq t\}$. Then, evidently, $s(t,g)+s(-t,g)\leq n(t,g)$,
and Jensen's formula gives
$$\int_0^r\frac{n(t,g)}{t}dt=\frac{1}{2\pi}\int_{-\pi}^\pi
\log|g(re^{i\theta})|d\theta+O(1)\leq \frac{2(b-a)}{\pi}(r+o(r)),
\quad r\to\infty.$$
The rest of the proof is the same as of Proposition~\ref{logans}.
\hfill$\Box$
\vspace{.1in}

The earliest results on the oscillation
of Fourier integrals with a spectral gap were obtained by M.G. Krein and
B.Ya. Levin in the 1940-s. The following result is contained
in \cite[Appendix II, Thm 5]{Levin}.
{\em Let $F\neq {\mathrm{const}}$ be a real function of bounded
variation
on the real line, and $dF$ has a spectral gap $(-a,a)$. Then}
\begin{equation}
\label{lev}
\liminf_{r\to\infty}\left\{ S(r,dF)-\frac{2a}{\pi}r\right\}>-\infty.
\end{equation}
This property neither follows from nor
implies (\ref{lowdens}).

In a footnote on p. 403 of \cite{Levin} Levin wrote:
``A similar, somewhat stronger result was obtained by M.G. Krein in the
theory of continuation of Hermitian-positive functions''.
Unfortunately, we were unable to find out what the precise formulation
of Krein's result was.

Application of a theorem of
Beurling--Malliavin as in the next section permits to
prove a version of Levin's theorem for functions satisfying
(\ref{logint}),
but we conjecture that (\ref{logint}) is not needed, that is
(\ref{log2}) holds for
arbitrary functions satisfying (\ref{exp}) with spectrum
on
$(-\infty,-a]\cup[a,\infty)$. Proposition~\ref{logans2}
supports this conjecture.

Recently Ostrovskii and Ulanovskii \cite{OU2}
extended and improved Levin's result
as follows:
{\em Let $dF(x)$ be a Borel measure satisfying
\begin{equation}
\label{ou}
\int\frac{|dF(x)|}{1+x^2}<\infty,
\end{equation}
where $|dF|$ stands for the variation,
and $dF$ has a
spectral gap $(-a,a)$. Then
$$\liminf_{r\to\infty}\left\{\int_1^r\left(\frac{1}{t^2}+\frac{1}{r^2}\right)
s(t,dF)dt-\frac{a}{\pi}\log r\right\}>0,$$
and}
$$\liminf_{r\to\infty}\left( S(r,dF)-\frac{2a}{\pi}r+3\log r\right)>0.$$
If $dF(x)=fdx,\; f\in L^1_{\mathrm{loc}}$, then
condition (\ref{ou})
is stronger than (\ref{logint}),
but weaker than the requirement of bounded variation in Levin's theorem.

These authors \cite{OU} also proved several
interesting results where the assumption about a spectral gap $(-a,a)$
is replaced by a weaker assumption that the Fourier transform of $f$
has an analytic continuation from the interval $(-a,a)$ to
a half-neighborhood of this interval in the complex plane.
However, in this result they characterize the oscillation of $f$ in terms of
the Beurling--Malliavin density
of sign changes, which is a sort of {\em upper density} rather than
lower density, see, for example, \cite[vol. II]{Koosis}.
\vspace{.1in}

The original results of Sturm in \cite{Sturm}
were about eigenfunctions of second order linear differential
operators $L$ on a finite interval; the case of trigonometric
polynomials corresponds to $L=d^2/dx^2$ on $[-\pi,\pi]$.
In 1916, Kellogg \cite{Kellogg} gave a rigorous proof of Sturm's claim
for
certain class of operators, whose inverses are defined
by totally positive symmetric kernels
on a finite interval $[a,b]$:
{\em Let $\phi_k$ be the $k$-th eigenfunction.
Then every linear combination
$$\sum_{k=m}^nc_k\phi_k\not\equiv 0, \quad n >m$$
has at least $m-1$ and at most $n-1$ sign changes on $[a,b]$.}

Our paper is based on a combination of two ideas; the first is the
proof of Sturm's theorem from \cite[III-184]{PS},
the second is similar to Sturm's own argument \cite[p. 430-433]{Sturm}
(compare \cite{Polya}).
We recall both proofs for the reader's
convenience.
\vspace{.1in}

1. Write the trigonometric polynomial (\ref{poly}) as
$$f(x)=h(x)+\overline{h}(x),\quad\mbox{where}\quad h(x)=\sum_{n\geq m}
c_ne^{inx},$$
then $h(x)=p(e^{ix})$ where $p$ is a
polynomial which has a root of multiplicity $m$ at zero.
By the Argument Principle, $p(z)$ makes at least $m$ turns
around  zero as $z$ describes the unit circle, so the curve
$\{ p(e^{ix}):0\leq x\leq2\pi\}$ intersects the imaginary axis
at least $2m$ times transversally. But $f(x)=2\Re h(x)$ changes
sign at each such intersection.
\vspace{.1in}

2. Use our trigonometric polynomial (\ref{poly}) as the initial
condition of the Cauchy Problem for the heat equation on
the unit circle. All coefficients will exponentially
decrease with time, and the lowest order term will have the slowest
rate of decrease. On the other hand,
as Sturm argued, the number of
sign changes of a temperature does not increase with time,
\cite{Sturm,Polya,Widder}.
So the number of sign changes of the initial condition is at least that
of the lowest degree term in its Fourier expansion.
\vspace{.1in}

In sections 3-4
we develop the first idea, and in sections 5-7 the second.

Other proofs of Sturm's theorem are given in \cite{PS}, problems
II-141, and VI-57.

To conclude this survey, we mention that Fourier Integral
first appears in Fourier's work on heat propagation \cite{Fourier1},
and that the study of sign changes was one of the main mathematical
interests of Fourier
during his whole career \cite{Fourier,Fourier1}.

\section{Application of the theorem of Beurling and Malliavin}

BMW are important because of the following
%
theorem of Beurling and Malliavin. {\em For every BMW $\omega$
and every $\eta>0$ there exists an entire function $g$ of exponential
type $\eta$,
such that
$g\exp\omega$ is bounded on the real line.}
The references are \cite{BM,Koosis2} and \cite[Vol. 2]{Koosis}.

Such function $g$ will be called an
$\eta$-multiplier. There is a lot of freedom in choosing a
multiplier, so we can ensure that $g$ has some additional properties.

First, there always exists a non-negative multiplier. Indeed, we can
replace $g$ by $g(z)\overline{g(\overline{z})}$.
A non-negative multiplier $g$ permits to reduce the proof of Theorem 1 to its
special case that $f\in L^1$. Indeed, let $f$ be a function satisfying
the conditions of Theorem 1. For arbitrary $\eta\in (0,a)$
we choose a non-negative $\eta$-multiplier $g$. Then $gf\in L^1,$ has the same
sequence of sign changes as $f$, and a spectral gap
$(-a+\eta,a-\eta)$. Applying Theorem 1 to $gf$ we obtain that
the sequence of sign changes of $f$ has lower density at least $(a-\eta)/\pi$,
for every $\eta\in (0,a)$. This implies (\ref{lowdens}).

We will use this observation in sections 5-7.

Second, there always exists a multiplier all of whose zeros are real.
(In fact, the multiplier constructed in the original
proof of the Beurling and Malliavin theorem has this property).
This we will use below in the proof of Proposition~\ref{111}.

Suppose that $f\in L^1$.
Then Fourier transform of $f$ is defined in the classical sense,
$$\hat{f}(t)=\int e^{-ixt}f(x)dx,$$
and $\hat{f}$ is a bounded function on the real line with the
property that $f(t)=0$ for $t\in (-a,a)$.
For $0<p<\infty$ we denote
$$\| h\|_p^*=\int\frac{|h(x)|^p}{1+x^2}dx,$$
and define the Hardy class $H^p$ as the set of all holomorphic
functions $h$ in the upper half-plane with the property that
$\| h(.+iy)\|_p^*$ is a bounded function of $y$ for $y>0$.
\begin{lemma}\label{hp} Let $f$ be a real function in $L^1$.
Then there exists a function $h$ in $H^{1/2}$
such that
\begin{equation}\label{main}
f(x)=h(x)+\overline{h}(x)\quad\mbox{a. e.},\quad\mbox{and}\quad h(iy)\to 0,\quad y\to+\infty,
\end{equation}
where $h(x)$ is the angular limit of $h$.
Furthermore,
\begin{equation}
\label{newh}
\| h\|_{1/2}^*\leq C_1\| f\|_1+C_2,
\end{equation}
where $C_1$ and $C_2$ are absolute constants.
If $f$ is an entire function then there is a unique entire function
$h\in H^{1/2}$ such that $(\ref{main})$ holds.
\end{lemma}

{\em Proof}. We define
\begin{equation}
\label{defh}
h(z)=\frac{1}{2\pi}\int_0^\infty e^{itz}\hat{f}(t)dt,\quad\Im z>0,
\end{equation}
which is evidently holomorphic in the upper half-plane.
Now we have for $\Im z>0$:
\begin{eqnarray*}
h(z)&=&\frac{1}{2\pi}\int_0^\infty e^{itz}\left\{\int e^{-its}f(s)ds\right\} dt\\
    &=&\frac{1}{2\pi}\int f(s)\left\{\int_0^\infty e^{it(z-s)}dt\right\} ds\\
    &=&\frac{i}{2\pi}\int f(s)\frac{ds}{z-s}.
\end{eqnarray*}
This representation shows that if $f$ is entire then
$h$ is entire: an analytic continuation of $h$ into the
lower half-plane can be obtained by deforming the path of integration in
the Cauchy integral.
Taking the real part, we obtain
$$2\Re h(x+iy)=\frac{y}{\pi}\int\frac{f(s)ds}{(x-s)^2+y^2}\;,
$$
so $2\Re h$ is the Poisson integral of $f$.
By Cauchy--Schwarz Inequality
\begin{equation}\label{cs}
\| \Re h\|^*_{1/2}\leq \sqrt{\pi}\| f\|_1^{1/2}.
\end{equation}
To prove that $h\in H^{1/2}$, we use the representation
of $\Im h$ as a Hilbert transform,
\begin{equation}
\label{hilb}
\Im h(x+iy)=\frac{1}{\pi}\int\left(\frac{1}{x-t}+\frac{t}{t^2+1}\right)
\Re h(t+iy)dt,
\end{equation}
and  Kolmogorov's inequality,
$$\displaystyle m(\lambda):=\int_{|\Im h(x+iy)|>\lambda}
\frac{dx}{1+x^2}\leq\frac{4}{\lambda}\int\frac{|\Re h(x+iy)|}{1+x^2}dx,$$
for each $\lambda>0$.
These can be found in \cite[v 1, p. 63]{Koosis}.
We have
\begin{eqnarray*}
&{}&\|\Im h(.+iy)\|_{1/2}^*=\int\frac{\sqrt{|\Im h(x+iy)|}}{1+x^2}dx\\
&=&-\int_0^\infty\lambda^{1/2}\,dm(\lambda)=\frac{1}{2}\int_0^\infty
\lambda^{-1/2}m(\lambda)d\lambda\\
&\leq &\pi+{2}
\|\Re h(.+iy)\|_1^*\int_1^\infty\lambda^{-3/2}d\lambda\leq C_1\|f\|_1
+C_2.
\end{eqnarray*}
Combined with (\ref{cs}), this implies (\ref{newh}).

Now we prove the uniqueness statement for the case that $f$ and $h$
are entire. Suppose that there are two representations
\begin{equation}
\label{rep}
f(x)=h_j(x)+\overline{h_j}(x),\quad j=1,2.
\end{equation}
The functions $w_+=h_1-h_2$ and $w_-=
\overline{h_2}-\overline{h_1}$ are analytic in the upper and lower
half-planes, respectively.
Subtracting one representation (\ref{rep}) from another, we obtain
$w_+(x)=w_-(x)$, so $w_+$ and $w_-$ are restrictions of a single entire
function $w$. By a theorem of Krein \cite[Ch. 16, Thm. 1]{LevB2}, $w$
is an efet of Cartwright class.
Now it follows from
(\ref{main}) that $w(iy)\to 0$ as $y\to\infty$
so $w=0$ and $h_1=h_2$.
\hfill$\Box$\vspace{.1in}

{\em Remark}. One can replace $H^{1/2}$ in Lemma~\ref{hp} by any $H^p$
with $p\in(0,1)$.

Now we restate the condition that $\hat{h}(t)=0$ for $t<a$ in terms of
$h$ itself.
\begin{lemma}\label{wp1} Let $h\in H^{1/2}$ be a function
represented by Fourier integral $(\ref{defh})$, where $\hat{f}$ is bounded,
and
$\hat{f}(t)=0$ for $t<a$.
Then $h$
satisfies
\begin{equation}\label{eee}
h(x+iy)=O(e^{-ay})\quad y\to\infty,
\end{equation}
uniformly with respect to $x$.
\end{lemma}

{\em Proof}.
$$
|h(x+iy)|\leq\frac{\|\hat{f}\|_\infty}{2\pi}\int_a^\infty e^{-sy}ds
\leq
\frac{e^{-ay}}{2\pi y}\| \hat{f}\|_\infty.$$
\hfill$\Box$\vspace{.1in}

We denote by $N$ the Nevanlinna class of functions of {\em bounded type}
in the upper half-plane. A holomorphic function $h$ in
the upper half-plane belongs to $N$ if $h$ is a ratio of bounded
holomorphic functions in the upper half-plane. We refer
to \cite{Nev,Privalov} for the theory of the class $N$.
Function $h$ from Lemma~\ref{hp} belongs to $N$ because $H^{p}\subset N$
for all $p>0$.
So we have the Nevanlinna representation
\begin{equation}\label{nev}
h(z)=e^{ia'z}B(z)e^{u(z)+iv(z)},
\end{equation}
where $a'$ is a real number,
$B$ a Blaschke product,
$u$ the Poisson integral of
$\log|h(x)|$, and $v$ the Hilbert transform of $u$ as in (\ref{hilb}).
In particular,
\begin{equation}
\label{int2}
J(u):=\int_{-\infty}^\infty\frac{|u(x)|}{1+x^2}dx<\infty.
\end{equation}
It is well-known that (\ref{nev}) implies
$$\limsup_{y\to+\infty}y^{-1}\log|h(iy)|= a',$$
so Lemma~\ref{wp1} gives $a'\geq a$.

To generalize Lemma~\ref{hp} to all functions satisfying (\ref{beu}) we
first recall the well-known fact:
\begin{lemma}\label{hurwitz}
Let $f$ be a function which satisfies $(\ref{exp})$,
and $g\in L^\infty(\R)$. If $f$ has a
spectral gap $(-a,a)$,
and $g$ is a function with spectrum
on $[-\eta,\eta],\;\eta< a,$ then $fg$ has a spectral gap
$(-a+\eta,a-\eta)$.
\end{lemma}
If $f$ also has bounded spectrum,
this follows from a theorem of Hurwitz,
\cite[Thm. 1.5.1]{Bieberbach}.
In the general case, the proof is the same; we include
it for the reader's convenience.

{\em Proof}. Let $F$ and $G$ be the Fourier transforms of $f$ and $g$.
Then $F$ is analytic in
$$\C\backslash\left\{ (-\infty,-a]\cup[a,\infty)\right\},$$
and $G$ is analytic in $\C\backslash[-\eta,\eta]$.
Let $\gamma$ be a simple closed curve going once counterclockwise
around the segment
$[-\eta,\eta]$, then
$$g(x)=\frac{1}{2\pi}\int_\gamma G(\zeta)e^{i\zeta x}d\zeta,$$
see, for example \cite{Levin,Bieberbach}.
We have
\begin{eqnarray*}
&&\int_0^\infty f(x)g(x)e^{-izx}dz
=\frac{1}{2\pi}\int_0^\infty
f(x)\int_\gamma G(\zeta)e^{i\zeta x}d\zeta e^{-izx}dx\\
&=&\frac{1}{2\pi}\int_\gamma G(\zeta)\int_0^\infty
f(x)e^{i(\zeta-z)x}dx\,d\zeta
=\frac{1}{2\pi}\int_\gamma G(\zeta) F(z-\zeta)d\zeta.
\end{eqnarray*}
This function is analytic in
$$\C\backslash\left\{(-\infty,-a+\eta]\cup[a-\eta,\infty)\right\}.$$
Similar computation for
$$-\int_{-\infty}^0 f(x)g(x)e^{-izx}dz$$
gives the same result. \hfill$\Box$
\vspace{.1in}


We state our conclusions as
\begin{proposition}\label{111} Let $f$ be a function satisfying the conditions
of Theorem 1. Then
$$f=h+\overline{h}\quad\mbox{a. e.,}$$
where $h$ is a function of bounded type in the upper half-plane,
having
representation $(\ref{nev})$ in which $a'\geq a$. If $f$ is an efet then
$h$ can be chosen in Cartwright's class.
\end{proposition}

{\em Proof}. Choose $\eta\in (0,a)$. Let $g$ be a Beurling-Malliavin multiplier
of exponential type $\eta$, real on the real line and having all zeros real.
Then $g$ is of bounded type in both upper and lower half-planes, and
\begin{equation}\label{bm}
\log|g(re^{i\theta})|=\eta r\sin\theta+o(r)\quad r\to\infty,
\end{equation}
uniformly with respect to $\theta$ for $|\theta|\in(\epsilon,\pi-\epsilon)$,
for every $\epsilon>0$.
Furthermore, $gf\in L^1$ by (\ref{beu}),
and $gf$ has a spectral gap $(-a+\eta,a-\eta)$ by Lemma~\ref{hurwitz}.
According to Lemma \ref{hp},
\begin{equation}\label{hp17}
gf(x)=h_1(x)+\overline{h_1}(x),
\end{equation}
where $h_1\in H^{1/2}$, so $h_1\in N$. Lemma \ref{wp1} implies that
\begin{equation}\label{as}
\log|h_1(re^{i\theta})|\leq(\eta-a)r\sin\theta+o(r)\quad r\to\infty,
\end{equation}
uniformly with respect to $\theta$.
Dividing (\ref{hp17}) by $g$ (which has no zeros outside the real axis),
we conclude that
(\ref{main}) holds with $h=h_1/g$ which evidently belongs to $N$.
Now (\ref{bm}) and (\ref{as}) show that
$$\log|h(re^{i\theta})|\leq -ar\sin\theta+o(r)\quad r\to\infty,$$
uniformly with respect to $\theta$,
which implies that $a'\geq a$ in (\ref{nev}).

If $f$ is an efet, let $b$ be its exponential type and
$F$ be its Fourier transform in the sense of Carleman.
Then $F$ is analytic in
$$\bC\backslash([-b,-a]\cup[a,b])$$
and $F(\infty)=0$. By the theorem on separation of singularities,
$F=F_1+F_2$, where $F_1$ is analytic
in $\bC\backslash[-b,-a]$, $F_2$ is analytic in $\bC\backslash[a,b]$,
and $F_j(\infty)=0,\; j=1,2$. This leads to the decomposition
\begin{equation}
\label{nnn}
f(x)=h^+(x)+h^-(x),
\end{equation}
where $h^{\pm}$ are efet with spectra on $[-b,-a]$ and $[a,b]$
respectively, so $h^{\pm}(iy)=O(\exp(-a|y|),\; y\to\pm\infty$.
Multiplying (\ref{nnn}) by $g$, and using the uniqueness statement
in Lemma \ref{hp} we obtain $h^+=h$ and $h^-=\overline{h}$, where
$h$ is a function of the class $N$ as above.
\hfill$\Box$
\vspace{.1in}

{\em Construction of Example 1}. We combine Logan's method
\cite[Thm 5.5.1]{Logan} with
the theorem of Beurling and Malliavin. Without loss of generality,
we may assume that $a=\pi-2\epsilon,$ and $b=\pi+2\epsilon$, where
$\epsilon>0$.
Let $g_1$ be a real entire function of zero exponential type,
satisfying (\ref{logint}), with only simple zeros,
and such that the zero set of
$g_1$ coincides with the set of integer points on the
intervals $[y_k,x_k]$:
$$g_1(n)=0, \; g'(n)\neq 0\quad\mbox{for}\quad n\in
\Z\cap\left(\cup_{k=1}^\infty [y_k,x_k]\right).$$
Such function $g_1$ can be easily constructed if the intervals
$[y_k,x_k]$ are not too long in comparison with $x_k$, for example,
if
$$\sum_{k=1}^n (x_k-y_k)\leq{x_n^\alpha}\quad\mbox{for some}\quad\alpha\in
(0,1).$$
One can obtain longer intervals, if desirable, whose size can
be characterized in terms of Beurling--Malliavin density \cite[vol. II]{Koosis}.
Let $g$ be an entire function of exponential type $\epsilon$, which is positive
on the real line and such that $|x|^2g(x)g_1(x)$ is bounded for
$x\in\R$. Such function $g$ exists by the Beurling--Malliavin theorem
(ii). Then
$$f_1(z)=g(z)g_1(z)\sin \pi z.$$
does not change sign on any of the intervals $[y_k,x_k]$,
and $\hat{f_1}$ has support on
$$[-\pi-\epsilon,-\pi+\epsilon]\cup[\pi-\epsilon,\pi+\epsilon].$$
Evidently, $f_1\in L^1$.
To destroy the multiple zeros of $f_1$ on the intervals $[y_k,x_k]$,
we put $f(z)=f_1(z+1/2)+f_1(z)$.
\hfill$\Box$

\section{Theorem 1 for real analytic functions}

To present the ideas unobscured by technical details,
we prove in this section Theorem~1 for real analytic functions $f$
whose real zeros are simple, so that the sign changes occur
exactly at the zeros of $f$.
The general case will be obtained from this special case in sections 5--7,
by a smoothing procedure.

We write, as in Proposition~\ref{111},
\begin{equation}\label{again}
f(x)=h(x)+\overline{h}(x),
\end{equation}
where $h$ has spectrum on $[a,\infty)$, and consider
the Nevanlinna representation (\ref{nev}).
Our assumptions about analyticity and simple zeros imply that
$v$ in (\ref{nev}) is piecewise continuous, the only jumps of
$-\pi$ occur exactly at the real zeros of $h$ (which are all simple).

Put
$$\phi(x)=\arg h(x):=a'x+\arg B(x)+v(x).$$
The Blaschke product
\begin{equation}\label{blaschke}
B(z)=\prod_{n\geq 1}\left(1-\frac{z}{z_n}\right)\left(1-
\frac{z}{\overline{z}_n}\right)^{-1},
\end{equation}
has a continuous argument because
zeros in the upper half-plane cannot accumulate to points on the real axis.
Furthermore, $\arg B$ is an increasing function, which is seen
by inspection of each factor of the product (\ref{blaschke}).

Let $\gamma$ be the curve in the $(x,y)$-plane consisting of the graph
of $\phi$ and vertical segments of length $\pi$ added at the points
of discontinuity of $v$. At each intersection of this curve with the set
\begin{equation}
\label{L}
L=\{ (x,y):x\in\R,\,y-\pi/2\in\pi\Z\},
\end{equation}
the number $h(x)$ is purely imaginary, that is $f(x)=0$ by (\ref{again}).

So we want to estimate from below the number of intersections
of $\gamma$ with $L$ over the intervals $[0,r]$.

We fix $\epsilon\in(0,1/2)$ and prove that on every interval
$[(1-\epsilon)x,x]$ with $x$ large enough there exists a point $x'$ such that
\begin{equation}\label{want}
\phi(x')\geq a'x'+v(x')>(a'-2\epsilon)x'.
\end{equation}
It will immediately follow from (\ref{want}) that the number of intersections
$\gamma\cap L$ has lower density at least $a'/\pi$. So it remains to
prove (\ref{want}).

We recall that $v$ is harmonically
conjugate to $u$, and that $u$ satisfies (\ref{int2}).
According to Kolmogorov's inequality \cite[v 1, p. 63]{Koosis}
$$\int_{|v(x)|>\lambda}\frac{dx}{1+x^2}\leq\frac{4}{\lambda}
\int_{-\infty}^\infty\frac{|u(x)|}{1+x^2}dx,$$
for each $\lambda>0$.

We
break $u$ into two parts with disjoint supports, $u=u_0+u_1$,
where the support of $u_0$
belongs to $[-r_0,r_0]$ for some $r_0>0$ and
$u_1$ satisfies
\begin{equation}\label{Ju1}
\int\frac{|u_1(x)|}{1+x^2}dx=
\int_{|x|>r_0}\frac{|u_1(x)|}{1+x^2}dx<\epsilon^2/8,
\end{equation}
which is possible in view of (\ref{int2}).
Let $v_j=\H u_j,\; j=0,1;$ where $\H$ stands for the Hilbert transform,
$$\H u(x)=\lim_{y\to 0+}\frac{1}{\pi}\left(\int\frac{x-t}{(x-t)^2+y^2}+
\frac{t}{t^2+1}\right)u(t)\,dt.$$
\begin{lemma}\label{lemma2}
$\displaystyle
|v_0(x)|\leq J(2r_0+r_0^{-1})/\pi\quad\mbox{for}\quad |x|>2r_0,$
where $J=J(u)$ is defined in $(\ref{int2})$.
\end{lemma}

{\em Proof}.
\begin{eqnarray*}
|v_0(x)|&\leq&\frac{1}{\pi}\left|\int_{-r_0}^{r_0}
\frac{u_0(t)}{x-t}dt\right|+\frac{1}{\pi}\left|\int_{-r_0}^{r_0}
\frac{tu_0(t)}{t^2+1}dt\right|\\
&\leq&\frac{1}{\pi}r_0^{-1}\int_{-r_0}^{r_0}|u_0(t)|dt+\frac{1}{\pi}r_0J\\
&\leq&\frac{1}{\pi}r_0^{-1}(1+r_0^2)J+\frac{1}{\pi}r_0J\\
&=&J(2r_0+r_0^{-1})/\pi.
\end{eqnarray*}
\hfill$\Box$

Now we prove that for every $x>2$ there exists
$$x'\in [(1-\epsilon)x,x],$$ such that
\begin{equation}\label{7}
v_1(x')>-\epsilon x.
\end{equation}
Suppose that this is not so. Then we apply Kolmogorov's inequality
to $v_1$ and $u_1$ with $\lambda=\epsilon x$, and (\ref{Ju1}):
$$\int_{(1-\epsilon)x}^x\frac{dt}{2t^2}<
\int_{(1-\epsilon)x}^x\frac{dt}{1+t^2}<\frac{4}{\epsilon x}
\int\frac{|u_1(x)|}{1+x^2}dx<\frac{\epsilon}{2x}.$$
Evaluating the integral on the left we conclude
$\epsilon/(1-\epsilon)<\epsilon,$ a contradiction.
This proves (our special case of) Theorem \ref{thm0}.
\vspace{.1in}

We state a more quantitative version of the result we just proved:
\begin{proposition}\label{quant}
Let $f$ be a function satisfying the conditions
of Theorem 1.
Suppose that $f$ is real analytic
and has only simple zeros on the real line. Write $f=h+\overline{h}$ as in
$(\ref{again})$, and let $h$ be represented by the formula $(\ref{nev})$,
with $J=J(u)$ as in $(\ref{int2})$.
Suppose that
$$\int_{|x|>r_0}\frac{|\log|h(x)||}{1+x^2}<\epsilon^2/8$$
for some $r_0>1$ and $\epsilon\in(0,1/2)$.
Then
$$s(r,f)\geq (a-\epsilon)r/\pi-J(2r_0+r_0^{-1})/\pi-1
\quad\mbox{for}\quad r>2r_0.$$
\end{proposition}
\hfill$\Box$

\section{Heating}
In this and the next two sections we assume that
$f\in L^1$ in Theorem 1. This does not restrict generality, as was explained in
the beginning of section 4.

If $f$ is not real analytic, or has multiple zeros
on the real line, we
``heat'' it.
This means that we
replace our $f$ by the
convolution\footnote{In the works
on heat equation this is called a Poisson integral.
We {\em don't do this}
to avoid confusion with the harmonic Poisson integral.}
with the heat kernel,
\begin{equation}
\label{conv}
f_t=K_t\ast f,\quad f_0=f,
\end{equation}
$$K_t(x)=\frac{1}{\sqrt{\pi t}}e^{-x^2/t}.$$
Evidently, $f_t$ are real analytic with respect to $x$ for all $t>0$.
All $$\hat{f}_t=\hat{K}_t\hat{f}=\exp(-s^2t/4)\hat{f}$$ have the same support
because
$\hat{K}_t$ never vanishes.

P\'olya \cite{Polya,Widder} proved that
$f_t$ has at most as many sign changes on
the real line as $f$ does. (This assertion was stated by Sturm
for the case of finite interval).
However, we cannot use this result\footnote{Probably it is
possible to derive what we need from P\'olya's result. However we think
it is useful
to give an independent proof of this generalization of
Sturm--P\'olya's theorem.}
because our functions have infinitely many sign changes,
and we have to control
their number on every interval $[0,r]$.
So we will prove the necessary generalization of P\'olya's theorem.

Our approach is closer to the original approach of Sturm rather than
that of P\'olya.

In this section we show that heating does not destroy
the conditions of Proposition \ref{quant}, and in the next two
sections we deal with the behavior of sign changes under heating,
and also with multiple roots which $f$ may have on the real line.
\begin{lemma}\label{logplus} Let $f\in L^1$ be a real
function with a spectral gap $(-a,a)$, and $f_t=K_t\ast f$.
Define $h_t$ by $(\ref{defh})$ using $f_t$ instead of $f$.

Then there exists $t_0>0$ such that
$\| h_t\|_1\leq\| h\|_1$, and $J(\log|h_t|)\leq C_1,$
for $t\in(0,t_0)$, and
where $C_1$ is independent of $t$.
Further, for every $\epsilon>0$ there exist $r_0>0$
such that for all $t\in(0,t_0)$ we have
\begin{equation}\label{17}
\int_{|x|\geq r_0} \frac{\left|\log|h_t(x)|\right|}{1+x^2}dx<\epsilon.
\end{equation}
\end{lemma}

We emphasize that $r_0$ and $C_1$ are independent on $t$.
They only depend of $h$ and $\epsilon$.

{\em Proof}.
First, of all,
$$\int|f_t(x)|dx=\int|K_t\ast f|(x)dx\leq\int(K_t\ast|f|)(x)dx=\int|f(x)|dx,$$
so $\| f_t\|_1\leq\| f\|_1$.
Using Lemma~\ref{hp} we obtain $\| h_t\|_{1/2}^*\leq C$, with $C$
independent of $t$. Thus
\begin{equation}\label{local}
\sqrt{|h_t(x)|}=k_t(x)(1+x^2),\quad\mbox{where}\quad \|k_t\|_1\leq C.
\end{equation}
We have
\begin{equation}\label{div}
\displaystyle
\begin{array}{rcl}
\log^+|h_t|&\leq&2\log^+|k_t|+2\log(1+x^2)\\
&\leq&2|k_t|+2\log(1+x^2)
\end{array}
\end{equation}
Let $u_t(x)=\log|h_t(x)|$ for real $x$ and $t\geq 0$.
Dividing (\ref{div}) by $1+x^2$,
integrating and using (\ref{local}) gives
\begin{equation}\label{6}
J(u_t^+)=
\int\frac{u^+_t(x)}{1+x^2}dx<C,
\end{equation}
where $C$ is independent of $t$.
Similarly we obtain from (\ref{div})
that
$$\int_{|x|\geq r_0}\frac{u^+_t(x)}{1+x^2}dx<\frac{2}{\sqrt{1+r_0}}
\left(\|k_t\|_1+
\int_0^\infty\frac{\log(1+x^2)}{(1+x^2)^{3/2}}dx\right)<\epsilon,$$
with some $r_0>1$ independent of $t$.

Property (\ref{6}) makes possible to extend $u^+_t$
to the upper half-plane
by Poisson's formula.
We continue to denote the extended function by $u^+_t$.
Notice that $h_t\in N$ for all $t$, and
$u^+_t(x+iy)-ay$ is a positive harmonic majorant
of $\log|h_t|$ in the upper half-plane.

Now we prove
\begin{equation}
\label{minus}
J(u_t^-)<C,
\end{equation}
with $C$ independent of $t$, and
\begin{equation}\label{34}
\int_{|x|\geq r_0}\frac{u^-_t(x)}{1+x^2}dx<\epsilon
\end{equation}
for some $r_0>0$.
Fix a point $z_0$ in the upper half-plane, such that $\delta=|h(z_0)|>0$.
As $h_t(z_0)\to h(z_0)$ as $t\to0$, we conclude that $h_t(z_0)>\delta/e$
when $t$ is small enough.
Let $b$ be the {\em true left end} of the support of $\hat{h}_t$.
It is important to notice that $b$ is {\em independent} of $t$,
because $\hat{h}_t=\hat{K}_t\hat{h}$.
Then
\begin{equation}\label{3}
u_t(z_0)-b\Im z_0\geq\log\delta-1>-\infty,
\end{equation}
when $t$ is small enough. Here we mean that $u_t$  is extended to a harmonic
function in the upper half-plane by the Poisson integral. Now (\ref{3})
implies (\ref{minus}). It remains to prove
(\ref{17}) for the negative part of $u_t$.
For psychological reasons it is better to work in the unit disc $\U$ instead
of the upper half-plane. The fractional-linear
transformation $T(z)=(z-i)/(z+i)$ maps the upper half-plane
onto $\U$, $T(\infty)=1$,
and we put
$\zeta_0=T(z_0),$ and
\begin{equation}\label{defw}
w_t=u_t\circ T^{-1}-b\Im T^{-1}.
\end{equation}
As a consequence of (\ref{3}) we have
\begin{equation}\label{4}
w_t(\zeta_0)\geq\log\delta-1>-\infty.
\end{equation}
The measure $dx/(1+x^2)$
on the real line corresponds to the measure $d\theta$ on
the unit circle $\T=\{ e^{i\theta}:\theta\in\R\}$.

It follows from (\ref{6})
that each $w_t$ is a difference of
positive harmonic functions in the unit disc, so it is the Poisson integral
of some charge $\mu_t$ of bounded variation on the unit circle.
The constant $b$ in (\ref{defw}) comes from the Nevanlinna representation
of $h_t$ similar to (\ref{nev}), and this constant
does not depend on $t$. So all charges $\mu_t$ have an atom of mass
exactly $-b$ at the point~$1$.

Let $\mu_t=\mu_t^+-\mu_t^-$ be the Jordan decompositions. Conditions
(\ref{4}) and (\ref{6}) imply that $\mu_t$ are of bounded total variation,
with a bound independent of $t$.
So we have weak convergence $\mu_t\to\mu_0,\; t\to 0$.
Let $\phi$ be a positive continuous function on the unit circle,
which is identically equal to $1$ in some neighborhood of the point $1$,
and at the same time
$$\left|b-\int_\T\phi|\mu_0|\right|<\epsilon/2,$$
where $|\mu_0|=\mu_0^++\mu_0^-$ is the variation of $\mu_0$.
Then there exists $t_0$ such that
\begin{equation}
\label{t0}
\left|b-\int_\T\phi|\mu_t|\right|<\epsilon,
\quad\mbox{for}\quad 0\leq t\leq t_0.
\end{equation}
When translated back to the real line from the unit circle, this implies
(\ref{34}).

\hfill$\Box$

\section{Preliminaries on temperatures}

Here we collect for the reader's convenience some
facts about convolutions (\ref{conv}) of real functions with the heat
kernel. We use the convenient notation\footnote{We apologize
for such abuse of the letter $u$, but the harmonic function
$u$ of sections
3-5 will not appear anymore until the end of section 7.}
\begin{equation}\label{tempu}
u(x,t)=f_t(x)
\end{equation}
and consider $u$ in the upper half-plane $\{ s=(x,t):t\geq 0,
x\in\R\}$.

The function $u$ in (\ref{tempu}) is a solution of the heat
equation in the open upper half-plane:
\begin{equation}\label{heateq}
4\frac{\partial u}{\partial t}=\frac{\partial^2 u}{\partial x^2}.
\end{equation}
Such functions are called {\em temperatures}. Formula (\ref{conv})
solves the initial value problem on an infinite rod (the $x$-axis)
with given initial temperature $f(x)$. A standard reference on the
subject is
\cite{Doob}.
Here is the precise statement about the boundary behavior of $u$
which is a slight generalization of
\cite[1.XVI.7]{Doob}:
\begin{lemma}\label{boundary}
Let $f$ be a real function from $L^1$.
Then for every $x\in\R$,
$$\liminf_{t\to 0}u(x,t)\geq\liminf_{\epsilon\to 0+}
\frac{1}{2\epsilon}\int_{x-\epsilon}^{x+\epsilon}f(t)dt.$$
\end{lemma}
This is a general property of positive symmetric kernels. Radial limits
can be replaced by non-tangential limits, and even by limits
from within parabolas tangent to the real line at $x$. It follows that
at every Lebesgue density point $x$ of $f$, the limit
$\lim_{t\to 0+}u(x,t)$ exists and equals $f(x)$.

Next lemma (due to L.~Nirenberg)
is called the Strong Minimum Principle \cite[1.XV.5]{Doob}
\begin{lemma}\label{strongmin}
Let $D$ be a bounded region in the horizontal strip $P=\{ s=(x,t):0<t<T\}$,
and $u$ a temperature in $D$. Suppose that
\begin{equation}\label{min}
\liminf_{ s\to \sigma}u(s)\geq 0,\quad\mbox{for all}\quad
\sigma\in \partial D\cap (P\cup(\R\times\{0\})).
\end{equation}
Then $u\geq 0$ in $D$, and if $u(s)=0$ for some point $s\in D$
then $u\equiv0$ in $D$.
\end{lemma}
We need an extension of the Minimum Principle, analogous to
the Phragm\'en--Lindel\"of Theorem in the theory of harmonic functions:
\begin{lemma}\label{phl}
Let $D$ be a region as in Lemma \ref{strongmin}, and $u$ a temperature in
$D$. Suppose that $u$ is bounded from below, and $(\ref{min})$ holds
for all but finitely many points $\sigma\in\partial D\cap(\R\times\{0\})$,
and the finite set of exceptional points belongs to the real axis.
Then the same conclusion as in Lemma \ref{strongmin} holds.
\end{lemma}

{\em Proof}. Let $x_1,x_2,\ldots,x_n$ be the exceptional points on
the real axis. Consider the auxiliary function
$$\displaystyle w(s)=\left\{\begin{array}{ll}
\displaystyle\sum_{k=1}^n\log^+\frac{1}{|s-x_k|},& s\in\R\times\{0\},\\
&\\
(K_t\ast w(.,0))(x),& s=(t,x),\; x\in\R, t>0.
\end{array}\right.
$$
Then $w$ is a positive temperature in $P$,
and $$w(s)\to+\infty\quad\mbox{as}\quad s\to x_k,\;\;s\in P,\;\;
1\leq k\leq n.$$
So, for every $\epsilon>0$, the function
$$u_\epsilon=u+\epsilon w$$
satisfies all conditions of Lemma \ref{strongmin}. So $u_\epsilon\geq 0$,
that is $u(z)\geq -\epsilon w(z)$. Letting $\epsilon\to 0$, we conclude
that $u\geq 0$. So $u$ satisfies the conditions of Lemma \ref{strongmin},
and the conclusions of Lemma \ref{strongmin} hold for $u$.
\hfill$\Box$\vspace{.1in}

\begin{lemma}\label{lemma6}
Let $u$ be a temperature in some region $D$ of $(x,t)$-plane.
Then multiple zeros of the functions $x\mapsto u(x,t)$
are isolated in $D$.
\end{lemma}

{\em Proof}. Suppose that $u$ has a non-isolated multiple zero,
Let $m\geq 2$ be the minimum of multiplicities of such zeros.
Then there exists an analytic germ $g(t)$ which gives
the position of such multiple zero for $t\in(t_0-\epsilon,
t_0+\epsilon)$ for some $t_0$ and $\epsilon>0$. So we have
$$u(x,t)=(x-g(t))^mv(x,t),$$
in a neighborhood of $(g(t_0),t_0)$. Here
$v$ is a real analytic
function
$$v(g(t_0),t_0)\neq 0.$$
We differentiate, and see that the lowest order term in
$\partial^2u/\partial x^2$ is $$m(m-1)(x-g(t))^{m-2}v(x,t),$$ while
all terms in $\partial u/\partial t$ are of order at least $m-1$.
So $u$ cannot satisfy the heat equation.
\hfill$\Box$\hspace{.1in}

\section{Heating, Part II}

In this section we complete the proof of Theorem 1.
Let $f\neq 0$ be a real function in $L^1(\R)$,
such that its Fourier transform $\hat{f}$ has a gap
$(-a,a)$. Choose an arbitrary $\epsilon>0$. Let $r_0$ be the number
defined in Lemma \ref{logplus}.
We will estimate the number of sign changes of $f$
on the interval $[0,r]$, where $r>2r_0$.
If $f$ has infinitely many sign changes on
$[0,r]$ then there is nothing to prove. So we assume that
the number of sign changes is finite on $[0,r]$.
A {\em zero place} of $f$ is defined as a maximal closed interval $I$,
(which may degenerate to a point)
such that $f|_I=0$ a.e. A zero place $I=[c,d]$ is called a
{\em place  of sign change} of $f$ if $(x-c)f(x)$ has constant
sign in a neighborhood of $I$. The complement of the union of
the places of sign changes consists of open intervals
which are called {\em intervals of constancy of sign}.
We write $I_1<I_2$ to mean that the
intervals $I_1$ and $I_2$ are disjoint and $I_2$ is on the right of
$I_1$.

Let $0<I_1<I_2<\ldots<I_n<r$ be the places of sign changes.
We assume that
\begin{equation}\label{assum}
 n\geq 2,\; f(0)<0,\; f(r)<0,
\end{equation}
and that $0$ and $r$ are
Lebesgue density points of $f$.
This assumption does not restrict generality.

Let $f_t=K_t\ast f$, and let $t_0$ be the number from
Lemma \ref{logplus}.
We are going to show, that for $t_0$ small enough,
the number of sign changes of $f_t$ on $[0,r]$ does not exceed that
of $f$ for $t\in(0,t_0)$.

Using Lemma \ref{boundary} and negativity of $f$ at its Lebesgue
points $0$ and $r$, we achieve that
\begin{equation}\label{inf}
\sup \{ f_t(x):x\in \{ 0,r\}, 0<t<T\}<0,
\end{equation}
by choosing $T\in(0,t_0)$ small enough.
We recall that $f_T$ is real analytic.
Using Lemma~\ref{lemma6}
we ensure that $f_T$ has only simple zeros.

{\em We are going to prove that}
\begin{equation}\label{claim}
s(r,f_T)\leq s(r,f).
\end{equation}

Assume first that $f$ is bounded in some neighborhood of the union
$\cup_{k=1}^nI_k$. As $f_T$ is real analytic, every place of sign change
of $f_T$ is one point.
We consider a maximal interval $\ell=(y_1,y_2)\subset[0,r]$ of sign
constancy of $f_T$, where $f_T$
is non-negative, but $y_1$ and $y_2$ are the places
of sign changes of $f_T$.
Define the strip $P=\{ s=(x,t):0<t<T\}$.
Denote, as in (\ref{tempu}),
$u(s)=u(x,t)=f_t(x)$.
Let $D$ be the connected component of the
set
$$\{ s\in P:u(s)>0\},\quad\mbox{such that}\quad\partial D\supset \ell.$$
Notice that $u(s)=0$ for $s\in\partial D\cap P$.
Then $D$ is bounded because it is contained in the
rectangle
$$\{ (s,t): 0<x<r,\, 0<t<T\}$$
in view of (\ref{inf}).
We claim that
\begin{equation}
\label{claim2}
\partial D\cap \{(x,t):t=T\}=\overline{\ell}=[y_1,y_2]\times\{T\}.
\end{equation}
Indeed, on those two intervals $\ell^-$ and $\ell^+$ of constant sign which
are adjacent to $\ell$, the sign is negative, so these
two intervals cannot intersect $\partial D$. If there is an
interval, say $\ell^*$, on the line $t=T$, which belongs to $\partial D$,
and $\ell^*\cap\ell=\emptyset$, we suppose,
for example that $\ell^*$ is on the same side of $\ell$ as $\ell^+$.
But then the component $D^+$ of the set
$\{ s\in P:u(s)<0\}$ which has $\ell^+$ on the boundary
has closure in the upper half-plane (being separated by $D$ from
the $x$-axis), and this contradicts Lemma \ref{strongmin}.
This proves our claim (\ref{claim2}).
\begin{center}\begin{picture}(0,0)%
\includegraphics{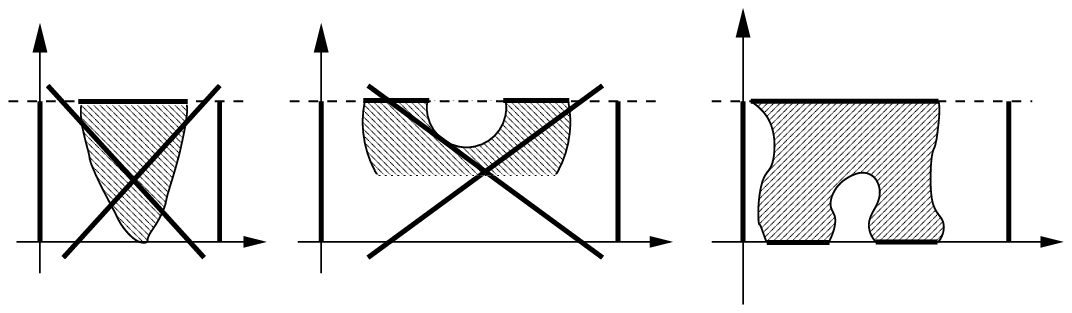}%
\end{picture}%
\setlength{\unitlength}{1973sp}%
\begingroup\makeatletter\ifx\SetFigFont\undefined%
\gdef\SetFigFont#1#2#3#4#5{%
  \reset@font\fontsize{#1}{#2pt}%
  \fontfamily{#3}\fontseries{#4}\fontshape{#5}%
  \selectfont}%
\fi\endgroup%
\begin{picture}(10169,3016)(129,-2462)
\put(6901,-256){\makebox(0,0)[lb]{\smash{\SetFigFont{12}{14.4}{\rmdefault}{\mddefault}{\updefault}
$T$
}}}
\put(2851,-256){\makebox(0,0)[lb]{\smash{\SetFigFont{12}{14.4}{\rmdefault}{\mddefault}{\updefault}
$T$
}}}
\put(151,-256){\makebox(0,0)[lb]{\smash{\SetFigFont{12}{14.4}{\rmdefault}{\mddefault}{\updefault}
$T$
}}}
\put(4374,-204){\makebox(0,0)[lb]{\smash{\SetFigFont{12}{14.4}{\rmdefault}{\mddefault}{\updefault}
$\ell_+$
}}}
\put(3826,-211){\makebox(0,0)[lb]{\smash{\SetFigFont{12}{14.4}{\rmdefault}{\mddefault}{\updefault}
$\ell^*$
}}}
\put(2101,-1936){\makebox(0,0)[lb]{\smash{\SetFigFont{12}{14.4}{\rmdefault}{\mddefault}{\updefault}
$r$
}}}
\put(5926,-1936){\makebox(0,0)[lb]{\smash{\SetFigFont{12}{14.4}{\rmdefault}{\mddefault}{\updefault}
$r$
}}}
\put(5109,-196){\makebox(0,0)[lb]{\smash{\SetFigFont{12}{14.4}{\rmdefault}{\mddefault}{\updefault}
$\ell$
}}}
\put(8101,-211){\makebox(0,0)[lb]{\smash{\SetFigFont{12}{14.4}{\rmdefault}{\mddefault}{\updefault}
$\ell$
}}}
\put(8153,-848){\makebox(0,0)[lb]{\smash{\SetFigFont{12}{14.4}{\rmdefault}{\mddefault}{\updefault}
$D$
}}}
\put(9676,-1936){\makebox(0,0)[lb]{\smash{\SetFigFont{12}{14.4}{\rmdefault}{\mddefault}{\updefault}
$r$
}}}
\put(1276,-211){\makebox(0,0)[lb]{\smash{\SetFigFont{12}{14.4}{\rmdefault}{\mddefault}{\updefault}
$\ell$
}}}
\put(8601,-2019){\makebox(0,0)[lb]{\smash{\SetFigFont{12}{14.4}{\rmdefault}{\mddefault}{\updefault}
$I_2$
}}}
\put(1276,-736){\makebox(0,0)[lb]{\smash{\SetFigFont{12}{14.4}{\rmdefault}{\mddefault}{\updefault}
$D$
}}}
\put(7550,-2033){\makebox(0,0)[lb]{\smash{\SetFigFont{12}{14.4}{\rmdefault}{\mddefault}{\updefault}
$I_1$
}}}
\put(3826,-886){\makebox(0,0)[lb]{\smash{\SetFigFont{12}{14.4}{\rmdefault}{\mddefault}{\updefault}
$D$
}}}
\end{picture}
\end{center}
If $\partial D\backslash\ell\subset P$, then $u\equiv 0$ in $D$ by
Lemma~\ref{strongmin},
and thus $u\equiv 0$ in the upper half-plane because
$u$ is real analytic.
If
$$\partial D\backslash\ell\subset P\cup(I_1\cup\ldots\cup I_n)\times\{0\},$$
we arrive at a contradiction in the similar way using Lemma \ref{phl}
instead of Lemma \ref{strongmin}. Here we used our temporary assumption
that $f$ was bounded in a neighborhood of $I_1\cup\ldots\cup I_n$.

The conclusion is that $\partial D$ intersects one of the
intervals $J$, a component of
the complement
 $$[0,r]\backslash\cup_{k=1}^nI_k.$$
But then $\partial D$ contains this interval
$J$ completely. This is because there is a neighborhood $U$ of $J$
such that $u(s)>0$ for $s\in U\cap P$, which follows from Lemma \ref{boundary}
combined with Lemma~\ref{strongmin}. Evidently, this $U$ cannot
intersect $\partial D\cap P$.

Thus $\partial D$ contains exactly one interval $\ell$ of
sign constancy of $f_T$ and at least one interval of
sign constancy of $f$. As different regions $D$ are evidently
disjoint, we conclude that $f$ has at least as many changes of signs
on $[0,r]$ as $f_T$. This proves (\ref{claim}).
\vspace{.1in}

It remains to get rid of the additional assumption that $f$ is bounded
in a neighborhood of $I_1,\ldots,I_n$. Let $U$ be a compact neighborhood
of these intervals in $\R$,
such that ${U}\cap\{0,r\}=\emptyset.$
For every positive integer $N$ we define
$$f^N(x)=\left\{\begin{array}{lll}
f(x)&\mbox{for}&x\not\in U,\\
f(x)&\mbox{if}&|f(x)|\leq N,\\
                                N&\mbox{if}&f(x)>N,\;x\in U\\
                                -N&\mbox{if}&f(x)<-N, \; x\in U.
\end{array}
\right.
$$
If $N$ is large enough, $f^N$ has the same number of sign changes
on $[0,r]$ as $f$.
Furthermore, $f^N\to f$ in $L^1$ as $N\to \infty$,
because $f(x)=f^N(x)$ for $x\not\in U$. As $U$ is disjoint from
the set $\{ 0, r\}$, the convergence $K_t\ast f^N\to f_t$
is uniform on $\{ 0, r\}\times[0,T]$, so for $N$ large enough
our functions $K_t\ast f^N$ are all strictly negative on this set.
So the previous proof applies to $K_T\ast f^N$,
and we conclude that $s(r,K_T\ast f^N)\leq s(r,f)$.
It remains to apply the observation that
$s(r,K_T\ast f^N)\to s(r,f_T),\; N\to\infty,$ pointwise, and thus
$s(r,f_T)\leq s(r,f)$.
So we proved (\ref{claim}) in full generality.
\vspace{.1in}

{\em Completion of the proof of Theorem 1}.
It remains to put the pieces together.
Let $f$ be a function satisfying the conditions of Theorem 1.
Assume wlog that $0$ is a Lebesgue point of $f$ and that $f(0)<0$.
Suppose, by contradiction, that for some $\eta\in (0,a/3)$ we have
$$\liminf_{x\to\infty}\frac{s(x,f)}{x}<\frac{a-3\eta}{\pi},$$
and let $x_k\to\infty$ be a sequence for which
\begin{equation}\label{llim}
s(x_k,f)<(a-3\eta)x_k/\pi.
\end{equation}
Apply the theorem of Beurling and Malliavin to find a multiplier
$g$ of type $\eta$, such that $g(x)\geq 0$ for $x\in\R$. Then
$gf\in L^1$ has the same sequence of sign changes as $f$, and a
spectral gap $(-a+\eta,a-\eta)$. We may assume that $x_k$ are
Lebesgue density points with $gf(x_k)<0$. For $t>0$, let
$(fg)_t=K_t\ast (fg)$ and let
$$(fg)_t=h_t+\overline{h_t}$$
be the decomposition which exists by Lemma \ref{hp}. Using Lemma \ref{logplus},
find
$r_0>0$ and $t_0>0$ such that
(\ref{17}) holds with $\epsilon=\eta^2/8$.
Choose $r=x_k>2r_0$ so that (\ref{llim}) is satisfied,
and
\begin{equation}\label{pizdec}
(a-2\eta)r/\pi-C_1(2r_0+r_0^{-1})/\pi-1>(a-3\eta)r/\pi,
\end{equation}
where $C_1$ is the upper bound for $J(\log|h_t|)$ from Lemma \ref{logplus}.
Then choose $t\in (0,t_0)$ so that
$(fg)_t$ has only simple zeros
on the real line (Lemma \ref{lemma6}), and the number
of these zeros on the interval $(0,r)$
is at most $s(r,f)=s(r,gf)$, which is guaranteed by (\ref{claim}).
Now by Proposition \ref{quant}, applied to $(gf)_t$, and (\ref{pizdec})
we have
$$s(r,f)=s(r,gf)\geq s(r,(gf)_t)>(a-3\eta)r/\pi,$$
where $r=x_k,$
which contradicts (\ref{llim}).
This proves the theorem.
\hfill$\Box$


\section{Limit sets of entire functions}

The theorem of Cartwright and Levinson
mentioned in sections 1 and 2 shows that
constructing an example of an efet whose indicator diagram is an interval
of the imaginary axis, and which does not have completely regular growth,
may be a non-trivial task. First such examples were constructed by
Redheffer and Roumieu \cite{R}, see also \cite{KR}. Their
purpose was to show that Titchmarsh's theorem on the
support of convolution fails for hyperfunctions with bounded support.
However, all these examples are still too regular for our purposes,
and we need the theory of limit sets, which generalizes
the theory of completely regular growth. It is due to
Azarin, Giner \cite{Azarin,AG,AzarinB},
H\"ormander and Sigurdsson \cite{HS}.
Here we collect the necessary facts from this theory.

Let $U^*$ be the set of all subharmonic functions in the plane
satisfying
$$\limsup_{|z|\to\infty}|z|^{-1}u(z)<\infty,$$
with induced topology from the space of Schwartz distributions $\D'(\C)$, and
$$U(\sigma)=\{ u\in U^*: u(0)=0,\; \sup_{z\in\C}|z|^{-1}u(z)\leq\sigma,\; \},$$
for $\sigma>0$. We recall that $\D'(\C)$ is a metric space.
We denote $U=\cup_{\sigma>0}U(\sigma).$

A one-parametric group $A$ of operators
$$(A_tu)(z)=t^{-1}u(tz),\quad t>0,$$
acts on $U^*$. The sets $U(\sigma)$ are $A$-invariant.

For a function $u\in U^*$ we define the {\em limit set} $\Fr\,[u]=\Fr_\infty[u]$
as the set of all limits
$$\lim_{n\to\infty}A_{t_n}u\quad\mbox{for}\quad t_n\to\infty.$$
Similarly, $\Fr_0[u]$ is defined for $u\in U$, using sequences $t_n\to 0$.
Each limit set $\Fr_\infty[u]$ or $\Fr_0[u]$
is a closed connected $A$-invariant subset
of $U(\sigma)$ for some $\sigma>0$.
If $f$ is an efet then $\log|f|\in U^*$, and we define the
{\em limit set of $f$} as $\Fr\,[\log|f|]$.
For every limit set $\Fr\,[u]$, the function
\begin{equation}
\label{indicator}
v(z)=\sup\{ w(z):w\in\Fr\,[u]\},
\end{equation}
is $A$-invariant and subharmonic. All such functions
have the form
\begin{equation}\label{ind}
v(re^{i\theta})=rh(\theta),\quad\mbox{where}\quad h''+h\geq 0,
\end{equation}
that is $h''+h$ is a non-negative measure. Functions $h$ with this property
are
called {\em trigonometrically convex}. The function $h$ defined by
(\ref{indicator}) and (\ref{ind}) is called the {\em indicator} of $u$.
If $f$ is an efet,
and $h$ the indicator of $\log|f|$ then $h$ coincides with the
classical Phragm\'en--Lindel\"of indicator of $f$.
The {\em indicator diagram}
is the closed convex set in the plane whose support function is $h$.

Criteria for a subset $\F\subset U$ to be a limit set
of some function $u\in U^*$
were found in \cite{AG} and \cite{HS}.
The following result is from \cite{AG} (see also \cite{AzarinB}).
\begin{proposition}\label{ag} Fix $\sigma>0$.
For a closed connected $A$-invariant subset
$\F\subset U(\sigma)$,
the following conditions are equivalent:
\newline
a) $\F=\Fr\,[u]$ for some $u\in U^*$,
\newline
b) $\F=\Fr\,[\log|f|]$ for some efet $f$, and
\newline
c) There exists a piecewise-continuous map
$$\R_{>0}\to U(\sigma),\quad t\mapsto v_t$$
with the properties
$$\dist(A_\tau v_t,\,v_{\tau t})\to 0,\quad t\to\infty,$$
and
$$\clos\{ v_t:t\in (t_0,\infty)\}=\F,\quad\forall t_0>0.$$
\end{proposition}
Here are some simple examples of limit sets derived from
Proposition~\ref{ag}.
\vspace{.1in}

\noindent
1. One-point limit set. Its only element has
to be of the form (\ref{ind}). This characterizes completely regular growth
in the sense of Levin--Pfluger.
\vspace{.1in}

\noindent
2. One periodic orbit. Let $u$ be a subharmonic function with
the property that $A_Tu=u$ for some $T\neq 1$. Then
$$\{ A_tu:1\leq t\leq T\}$$
is a limit set. One can show that in this case the indicator diagram
cannot be a non-degenerate interval of the imaginary axis,
so this type of functions
is not appropriate for our purposes.
\vspace{.1in}

\noindent
3. The closure of a single orbit,
$$\{ A_tu:0<t<\infty\}\cup\Fr_0[u]\cup\Fr_\infty[u],
\quad\mbox{where}\quad u\in U(\sigma)$$
is a limit set if and only if
$$\Fr_0[u]\cap\Fr_\infty[u]\neq\emptyset.$$
Again, in this case the indicator diagram cannot be
a non-degenerate interval of the imaginary axis.
\vspace{.1in}

\noindent
3. An interval. If $u_0$ and $u_1$ are two $A$-invariant
functions in $U$ then the set
$$\{ tu_0+(1-t)u_1:0\leq t\leq 1\}$$
is a limit set.
\vspace{.1in}

Examples in \cite{KR} are of this sort.
The efet constructed in \cite{KR} has indicator diagram $[-ib,ib]$
and the {\em lower} density of zeros is strictly less than $b/\pi$.
We need an example of efet with the indicator diagram $[-ib,ib]$
and the {\em upper} density of positive zeros strictly greater that $b/\pi$.
To achieve this we
combine the last two examples.
\begin{lemma}\label{Az}
Let $u$ be a function in $U$ with the properties
$$\Fr_0[u]=\{ u_0\}\quad\mbox{and}\quad\Fr_\infty[u]=\{ u_1\}.$$
Then
\begin{equation}
\label{ef}
\F=\{ A_tu:0<t<\infty\}\cup\{tu_0+(1-t)u_1:0\leq t\leq 1\}
\end{equation}
is a limit set.
\end{lemma}

{\em Proof}.
This easily follows from the general criterion in Proposition~\ref{ag}.
Fix a sequence of positive numbers with the property $r_{k+1}/r_k\to\infty,\;
k\to\infty$.

If $k=n^2$ for a positive integer $n$, we set $s_k=\sqrt{r_kr_{k+1}}$, and
$$v_t=A_{t/s_k}u,\quad r_k\leq t< r_{k+1}.$$
If $k=n^2+j,$ where $1\leq j\leq 2n$, we define
$$v_t=(j/2n)u_0+((2n-j)/2n)u_1,\quad r_k\leq t<r_{k+1}.$$
Then it is easy to verify that $v_t$ satisfies condition c) of Proposition~\ref{ag}
with $\F$ as in (\ref{ef}).
\hfill$\Box$
\vspace{.1in}

Now we describe the relation between the limit set and the distribution
of zeros of an efet.
Consider the set of all Borel measures in $\C$ (non-negative and
such that the measure
of every compact set is finite).
The analog of operators $A_t$ for measures
is
$$(B_t\mu)(E)=t^{-1}\mu(tE),\quad\mbox{for Borel sets}\quad E\subset\C.$$
Laplace operator $(2\pi)^{-1}\Delta$ splits $A_t$ and $B_t$:
\begin{equation}
\label{split}
\Delta A_t=B_t\Delta.
\end{equation}
We denote by $V^*$ the set of all measures $\mu$,
which satisfy
$$\limsup_{r\to\infty}r^{-1}\mu(D(r))<\infty,$$
where $D(r)=\{ z\in\C:|z|\leq r\},\; r\geq 0$.
We also define the subsets
$$V(\sigma)=\{\mu\in V^*:r^{-1}\mu(D(r))\leq\sigma,\; 0<r<\infty\},\quad
\sigma>0,$$
and $V=\cup_{\sigma>0}V(\sigma).$
Laplace operator is continuous in $U$ and sends
$U$ to $V$ (however, this map is not surjective,
and the image of $U(\sigma)$ is not equal to $V(\sigma')$
for any $\sigma'>0$).
Given a measure $\mu\in V^*$, we define the limit set $\Fr\,[\mu]$
as the set of all limits in $D'(\C)$
$$\lim_{n\to\infty}B_{t_n}\mu\quad\mbox{for}\quad t_n\to\infty.$$
It follows from (\ref{split}) that for every
$u\in U^*$ we have
\begin{equation}
\label{split2}
(2\pi)^{-1}\Delta\left(\Fr\,[u]\right)=\Fr\,\left[
(2\pi)^{-1}\Delta u\right].
\end{equation}
If $f$ is entire then $(2\pi)^{-1}\Delta\log|f|$ is the counting measure of zeros
of $f$. So the asymptotic distribution of zeros is reflected
in the Riesz measures of the elements of the limit set.
Let us make this more precise.
Two measures in $U^*$ are called equivalent if
$$B_t(\mu_1-\mu_2)\to 0\quad\mbox{as}\quad t\to\infty.$$
This implies $\Fr\,[\mu_1]=\Fr\,[\mu_2]$.
Let $T:\C\to\C$ be a map with the property
\begin{equation}
\label{prop}
T(z)-z=o(z),\quad z\to\infty.
\end{equation}
We recall that push-forward $T_*\mu$ of a measure by $T$ is defined by
$(T_*\mu)(E)=\mu(T^{-1}(E)).$
If $\mu\in V^*$, and a map $T$ satisfies (\ref{prop}), then $T_*\mu$ is
equivalent to $\mu$.
For each $\mu\in V^*$ one can construct
a map $T$ with the property (\ref{prop})
such that $T_*\mu$ is a counting measure of a divisor in $\C$.
This explains the implication a)$\rightarrow$b) in Proposition~\ref{ag}.
\begin{lemma}\label{Az2} Let $\mu$ be a measure in $V^*$.
Suppose that all measures in $\Fr\,[\mu]$ are supported on the real line
and have the form $d(x)dx$ where $d(x)<1$.
Then there exists a measure $\mu'$, which is equivalent to $\mu$
and which is supported on the integers, and $\mu'(n)\in\{0,1\}$ for
each integer $n$.
\end{lemma}

{\em Proof}. First we project our measure $\mu$ on the real line
by the map
$$T(re^{i\theta})=\left\{\begin{array}{ll}r,&|\theta|<\pi/2,\\
                                          -r,&|\theta-\pi|\leq\pi/2.
\end{array}\right.$$
This map does not satisfy
(\ref{prop}) but it is easy to see that
$\mu^{\prime\prime}=T_*\mu\sim \mu$ for measures $\mu$ satisfying the conditions
of Lemma~\ref{Az2}

Second,  let $F$ be the distribution function of $\mu^{\prime\prime}$,
that is $\mu^{\prime\prime}=dF$ and
$F(0)=0$. Then we set $F_1(x)=F([x])$, where $[.]$ stands
for the integer part, and put
$\mu'=dF_1$.
It is clear that the jumps of $F_1$ are at most $1$, and they occur
only at integers.
\hfill$\Box$

\section{Example of a hyperfunction}

Here we construct Example 2 assuming, without loss of generality,
that
$a+b=2\pi$.
We begin with a smooth negative function $u=ku_0$ with support on $[0,2]$,
for example, we can take
$$\displaystyle u(x)=\left\{\begin{array}{ll}
       -k(1-(x-1)^2)^2,& |x-1|\leq 1,\\
     \\
       0,& |x-1|>1,\end{array}\right.$$
where $k>0$ is a parameter to be specified later.
Then we extend $u$ to $\C\backslash\R$ by Poisson's integral.
The resulting function $u$ is a delta-subharmonic function in $\C$,
whose Riesz charge is supported on $\R$ and has the form
$dQ(x)=q(x)dx$, where $q$ is a smooth function.
So we have,
$$u(z)=\int\log\left|1-\frac{z}{t}\right|dQ(t)=\int\log\left|
1-\frac{z}{t}\right|q(t)dt.$$
We notice that $u|_\R$ is the Hilbert transform of $Q$.
So for the function $u$ as above we can explicitly compute $Q$ and $q$:
$$Q(x+1)=k(x^2-1)^2\log\left|\frac{x+1}{x-1}\right|-2kx^3+\frac{10}{3}kx,$$
and
$$q(x+1)=4kx(x^2-1)\log\left|\frac{x+1}{x-1}\right|-8kx^2+\frac{16}{3}k.$$
We put
\begin{equation}
\label{m}
-m=\min_{x\in\R}q(x)<0,
\end{equation}
and
\begin{equation}
\label{eta}
\eta=\max_{x\geq 0}\frac{Q(x)}{x}>0.
\end{equation}
Now we choose and fix $k$ so small that
\begin{equation}
\label{mpluseta}
m+\max_{x\in\R}q(x)<1.
\end{equation}
We define
\begin{equation}
\label{Q1}
Q_1(x)=Q(x)+mx,\quad\mbox{so that}\quad q_1=Q_1^\prime=q+m\geq 0,
\end{equation}
in view of (\ref{m}), and thus the function
\begin{equation}
\label{u1}
u_1(z)=\int\left(\log\left|1-\frac{z}{t}\right|+\Re\left(\frac{z}{t}\right)\right)dQ_1(t)
=u(z)+\pi m|\Im(z)|,
\end{equation}
is subharmonic in $\C$ and belongs to the class $U$ defined in the
previous section.
We have
\begin{equation}
\label{fru1}
\Fr_\infty[u_1]=\{ \pi m|\Im (.)|\}\quad\mbox{and}\quad
\Fr_0[u_1]=\{\pi m'|\Im(.)|\},
\end{equation}
where $\Im(.)$ is the function $z\mapsto\Im(z)$, and
\begin{equation}
\label{nn}
m'=m+q(0)<m.
\end{equation}
The first formula in (\ref{fru1}) follows from $u(z)\to 0$ as $z\to\infty$,
while the second one and (\ref{nn}) follow from
$q(0)=Q'(0)=(\partial u/\partial y)(0)$.
Now by Lemma~\ref{Az}, the set
$$\F:=\{ A_tu_1:t\in\R\}\cup\{ t|\Im(.)|:\pi m'\leq t\leq \pi m\}\subset U$$
is a limit set of an efet.
Evidently,
\begin{equation}
\label{indic}
\sup\{ w(z):w\in\F\}=\pi m|\Im(z)|.
\end{equation}
Let $g$ be an entire function of exponential type $m$,
such that
$$\Fr\,[\log |g|]=\F.$$
According to (\ref{indic}), the indicator diagram of $g$ is the
interval $[-\pi mi,\pi mi]$. In other words, Fourier transform of
$g$ is a hyperfunction supported on $[-\pi m,\pi m]$, \cite[v.2,
Thm.15.1.5]{Hormander}

In addition, we require that all zeros of $g$, except $o(r)$ of them
be simple and located at integers,
which is possible by Lemma~\ref{Az2} because the Riesz measures of
all elements of $\F$ are concentrated on the real line, and their densities
do not exceed $1$ in view of (\ref{mpluseta}).
The upper density of zeros of $g$
on the positive ray is
\begin{equation}\label{updens}
\max_{x\geq 0}Q_1(x)/x=m+\eta.
\end{equation}
Indeed, the limit set $\F=\Fr\,[\log|g|]$ contains $u$. This means that
there is a sequence $t_k\to\infty$ such that $B_{t_k}\Delta\log |g|\to\Delta u$;
this follows from (\ref{split2}).
Suppose that the maximum in (\ref{updens}) is attained at a point $x^*>0$.
Put $r_k=t_kx^*$ and let $n(r)$ be the number of zeros of $g$ on the interval
$[0,r]$. Then $n(r_k)/t_k=Q(x^*)+o(1),\; k\to\infty$,
and thus $n(r)/r=Q(x^*)/x^*+o(1)=m+\eta+o(1),\; r\to\infty$.

Finally we set
$$f(z)=g(z)\sin\pi z.$$
Then Fourier transform of $f$
$f$ is a hyperfunction
supported on
$$\pi[-1-m,\,-1+m]\cup\pi[1-m,\,1+m],$$
while the sign changes occur only at those integers which are
not zeros of $g$, that is the lower density of sign changes is
at most $1-m-\eta<1-m$ in view of (\ref{updens}).
\hfill$\Box$

{\em Purdue University

West Lafayette IN 47907}

{\em eremenko@math.purdue.edu

dmitry@math.purdue.edu}
\end{document}